\documentclass[draft,reqno]{amsart}
\usepackage{amsfonts}%
\usepackage{mathrsfs}

 \usepackage{cite}
\usepackage[dvips]{color}
\allowdisplaybreaks
\date{\today}

\theoremstyle{plain}
\newtheorem{thm}{Theorem}[section]

\newtheorem{lem}[thm]{Lemma}
\newtheorem{prop}[thm]{Proposition}

\theoremstyle{definition}

\theoremstyle{remark}
\newtheorem{rem}{Remark}[section]
\numberwithin{equation}{section}

\renewcommand{\d}{{\mathbf J}}
\renewcommand{\u}{{\mathbf u}}
\renewcommand{\v}{{\mathbf{v}}}
\renewcommand{\t}{{\mathbf{B}}}
\renewcommand{\H}{\mathbf{H}}

\newcommand{\w}{{\mathbf w}}
\newcommand{\R}{{\mathbb R}}

\newcommand{\dv}{{\rm div }}
\newcommand{\cu}{{\rm curl\, }}
\newcommand{\E}{{\mathcal E}}

\begin{document}

\title[Low Mach number limit for full magnetohydrodynamic equations ]
{Low Mach number limit for the multi-dimensional Full
magnetohydrodynamic equations}

\author{Song Jiang}
\address{LCP, Institute of Applied Physics and Computational Mathematics, P..O.
 Box 8009, Beijing 100088, P.R. China}
 \email{jiang@iapcm.ac.cn}

\author{Qiangchang Ju}
\address{Institute of Applied Physics and Computational Mathematics, P.O.
 Box 8009-28, Beijing 100088, P.R. China}
 \email{qiangchang\_ju@yahoo.com}

 \author[Fucai Li]{Fucai Li}
\address{Department of Mathematics, Nanjing University, Nanjing
 210093, P.R. China}
 \email{fli@nju.edu.cn}



\begin{abstract}
The low Mach number limit for the multi-dimensional full
magnetohydrodynamic equations, in which the effect of thermal conduction
is taken into account, is rigorously justified in the framework of classical
solutions with small density and temperature variations.
Moreover, we show that for sufficiently small Mach number,
the compressible magnetohydrodynamic equations admit a smooth
solution on the time interval where the  smooth solution of the
incompressible magnetohydrodynamic equations exists. In addition,
the low Mach number limit for the ideal magnetohydrodynamic equations with
small entropy variation is also investigated. The convergence rates are
obtained in both cases.
\end{abstract}

\maketitle

\section{Introduction}
The magnetohydrodynamic (MHD) equations govern  the  motion  of
compressible quasi-neutrally ionized fluids under the influence of
electromagnetic fields. The full three-dimensional compressible MHD
equations read as (see, e.g., \cite{Hu,JT64,KL,LL})
\begin{align}
&\partial_t\rho +\dv(\rho\u)=0, \label{naa} \\
&\partial_t(\rho\u)+\dv\left(\rho\u\otimes\u\right)+ {\nabla
p}
  =\frac{1}{4\pi}(\nabla \times \H)\times \H+\dv\Psi, \label{nab} \\
&\partial_t\H-\nabla\times(\u\times\H)=-\nabla\times(\nu\nabla\times\H),\quad
\dv\H=0, \label{nac}\\
&\partial_t\E+\dv\left(\u(\E'+p)\right)
=\frac{1}{4\pi}\dv((\u\times\H)\times\H)\nonumber\\
& \qquad \qquad \qquad\qquad\qquad\, \, +\dv\Big(\frac{\nu}{4\pi}
\H\times(\nabla\times\H)+ \u\Psi+\kappa\nabla\theta\Big).
 \label{nad}
\end{align}
Here  $x\in \Omega$, and $\Omega$ is assumed to be the whole
$\mathbb{R}^3$ or the  torus $\mathbb{T}^3$. The unknowns $\rho $
denotes the density, $\u=(u_1,u_2,u_3)\in \R^3$ the velocity, $\H=(H_1,H_2,H_3)\in \R^3$ the
magnetic field, and $\theta$ the temperature, respectively; $\Psi$
is the viscous stress tensor given by
\begin{equation*}
\Psi=2\mu \mathbb{D}(\u)+\lambda\dv\u \;\mathbf{I}_3
\end{equation*}
with $\mathbb{D}(\u)=(\nabla\u+\nabla\u^\top)/2$,
 $\mathbf{I}_3$ being the $3\times 3$ identity matrix,
and $\nabla \u^\top$   the transpose of the matrix
$\nabla \u$; $\E$ is the total energy given by $\E=\E'+|\H|^2/({8\pi})$ and
$\E'=\rho\left(e+|\u|^2/2 \right)$ with $e$ being the internal
energy, $\rho|\u|^2/2$ the kinetic energy, and
$|\H|^2/({8\pi})$ the magnetic energy. The viscosity
coefficients $\lambda$ and $\mu$ of the flow satisfy
$2\mu+3\lambda>0$ and $\mu>0$; $\nu>0$ is the magnetic diffusion
coefficient of the magnetic field,   and $\kappa>0$ is the heat
conductivity. For simplicity, we assume that $\mu,\lambda,\nu$ and
$\kappa$ are constants. The equations of state $p=p(\rho,\theta)$ and
$e=e(\rho,\theta)$ relate the pressure $p$ and the internal energy
$e$ to the density $\rho$ and the temperature $\theta$ of the flow.

The MHD equations have attracted a lot of attention of physicists and
mathematicians because of its physical importance, complexity, rich
phenomena, and mathematical challenges,
see, for example, \cite{BT02,CW02,CW03,Hu,DF06,FJN,FS,HT,LL,HW2,ZJX} and the
references cited therein. One of the important topics on the
equations \eqref{naa}--\eqref{nad} is to study its low Mach number limit.
For the isentropic MHD equations, the low Mach number limit has been rigorously
proved in \cite{KM,HW3,JJL1,JJL2}. Nevertheless, it is more significant and
 difficult to study the limit for the non-isentropic models from both physical and
 mathematical points of view.

The main purpose of this paper is to  present the rigorous
justification of the low Mach number limit  for the full MHD equations
\eqref{naa}-\eqref{nad} in the framework of classical solutions.

Now, we rewrite the energy equation (\ref{nad}) in the form of the internal energy.
Multiplying  \eqref{nab} by $\u$ and \eqref{nac} by $\H/({4\pi})$,
and summing them together, we obtain
\begin{align}\label{naaz}
&\frac{d}{dt}\big(\frac{1}{2}\rho|\u|^2+\frac{1}{8\pi}|\H|^2\big)
+\frac{1}{2}\dv\Big(\rho|\u|^2\u\Big)+\nabla p\cdot\u \nonumber\\
&\quad  =\dv\Psi\cdot
\u+\frac{1}{4\pi}(\nabla\times\H)\times\H\cdot\u
+\frac{1}{4\pi}\nabla\times(\u\times\H)\cdot\H\nonumber\\
& \, \, \qquad -\frac{\nu}{4\pi}\nabla\times(\nabla\times\H)\cdot\H.
\end{align}
Using  the identities
\begin{equation}
 \dv(\H\times(\nabla\times\H))  =|\nabla\times\H|^2-\nabla\times(\nabla\times\H)\cdot\H
 \end{equation}
 and
 \begin{equation}
\dv((\u\times\H)\times\H)   =(\nabla\times\H)\times\H\cdot\u+\nabla\times(\u\times\H)\cdot\H, \label{nae}
\end{equation}
and subtracting \eqref{naaz} from \eqref{nad}, we obtain the internal
energy equation
\begin{equation}\label{nagg}
\partial_t (\rho e)+\dv(\rho\u e)+(\dv\u)p=\frac{\nu}{4\pi}|\nabla\times\H|^2+\Psi:\nabla\u+\kappa \Delta \theta,
\end{equation}
where $\Psi:\nabla\u$ denotes the scalar product of two matrices:
\begin{equation*}
\Psi:\nabla\u=\sum^3_{i,j=1}\frac{\mu}{2}\left(\frac{\partial
u_i}{\partial x_j} +\frac{\partial u_j}{\partial
x_i}\right)^2+\lambda|\dv\u|^2=
2\mu|\mathbb{D}(\u^\epsilon)|^2+\lambda(\mbox{tr}\mathbb{D}(\u^\epsilon))^2.
\end{equation*}

In this paper, we shall focus our study on the ionized fluid obeying the perfect gas relations
\begin{align}
p=\mathfrak{R}\rho \theta,\quad e=c_V\theta,\label{pg}
\end{align}
where the constants $\mathfrak{R}, c_V\!>\!0$ are the gas constant and
the heat capacity at constant volume, respectively. We point out here that our analysis
below can be applied to more general equations of state for $p$ and $e$
by employing minor modifications in arguments.

To study the low Mach number limit of the system \eqref{naa}--\eqref{nac}
 and \eqref{nagg}, we use its appropriate dimensionless form as
follows (see the Appendix for the details)
\begin{align}
   &\partial_t\rho +\dv(\rho\u)=0, \label{naaa} \\
&\rho(\partial_t\u+\u\cdot\nabla\u)+\frac{\nabla (\rho\theta)}{\epsilon^2}
  =(\nabla \times \H)\times \H+\dv\Psi, \label{nabb} \\
&\partial_t\H-\nabla\times(\u\times\H)=-\nabla\times(\nu\nabla\times\H),\quad
\dv\H=0,\label{nacc} \\
&\rho(\partial_t
\theta+\u\cdot\nabla\theta)+(\gamma-1)\rho\theta\dv\u=\epsilon^2\nu|\nabla\times\H|^2
+\epsilon^2\Psi:\nabla\u+\kappa \Delta \theta,
  \label{nadd}
\end{align}
where $\epsilon= M$ is the  Mach number and the
coefficients $\mu,\lambda,\nu$ and $\kappa$ are the scaled
parameters. $\gamma=1+\mathfrak{R}/c_V$ is the ratio of specific
heats. Note that we have used the same notations and assumed that
the coefficients $\mu,\lambda,\nu$ and $\kappa$ are independent of
$\epsilon$ for simplicity. Also, we have ignored the Cowling number in the
equations \eqref{naaa}--\eqref{nadd}, since it does not create
any mathematical difficulties in our analysis.

 We shall study the limit as $\epsilon\to 0$ of the solutions to
 \eqref{naaa}--\eqref{nadd}.
We further restrict ourselves to the small density and temperature
variations, i.e.
\begin{align}
\rho=1+\epsilon q,\quad \theta=1+\epsilon \phi.\label{st1}
\end{align}

We first give a formal analysis. Putting \eqref{st1} and \eqref{pg}
into the system \eqref{naaa}--\eqref{nadd}, and using the identities
\begin{gather}
\cu\cu {\H}   =  \nabla\,\dv {\H}-\Delta {\H},\nonumber\\
\nabla(|{\H}|^2)   =2{\H}\cdot \nabla
{\H}+2{\H}\times \cu {\H},\label{idena}\\
 \cu( { \u}\times{\H})   =
 {\u} (\dv {\H})  - {\H}  (\dv{\u})
 + {\H}\cdot \nabla{\u} - {\u}\cdot \nabla{\H},\label{idenb}
\end{gather}
  then we can rewrite \eqref{naaa}--\eqref{nadd} as
\begin{align}
&\partial_t q^\epsilon +\u^\epsilon\cdot\nabla q^\epsilon
+\frac{1}{\epsilon}(1+\epsilon q^\epsilon)\dv\u^\epsilon=0, \label{na1} \\
&(1+\epsilon q^\epsilon)(\partial_t\u^\epsilon+\u^\epsilon\cdot\nabla
\u^\epsilon)+\frac{1}{\epsilon}\big[ (1+\epsilon q^\epsilon)
\nabla \phi^\epsilon+(1+\epsilon \phi^\epsilon)\nabla q^\epsilon\big]\nonumber\\
&\qquad\qquad-\H^\epsilon\cdot\nabla\H^\epsilon+\frac{1}{2}\nabla(|\H^\epsilon|^2)
=2\mu\dv(\mathbb{D}(\u^\epsilon))+\lambda\nabla(\mbox{tr}\mathbb{D}(\u^\epsilon)), \label{na2} \\
&\partial_t\H^\epsilon+\u^\epsilon\cdot\nabla\H^\epsilon +\dv
\u^\epsilon\H^\epsilon-\H^\epsilon\cdot\nabla
\u^\epsilon=\nu\Delta\H^\epsilon,\quad \dv\H^\epsilon=0,
\label{na4}\\
&(1+\epsilon q^\epsilon)(\partial_t
\phi^\epsilon+\u^\epsilon\cdot\nabla
\phi^\epsilon)+\frac{\gamma-1}{\epsilon}(1+\epsilon q^\epsilon)
(1+\epsilon \phi^\epsilon)\dv \u^\epsilon\nonumber\\
&\qquad\qquad=\kappa\Delta
\phi^\epsilon+\epsilon(2\mu|\mathbb{D}(\u^\epsilon)|^2
+\lambda(\mbox{tr}\mathbb{D}(\u^\epsilon))^2)
+\nu\epsilon|\nabla\times \H^\epsilon|^2.   \label{na3}
\end{align}
Here we have added the superscript $\epsilon$ on the unknowns to stress the
dependence of the parameter $\epsilon$. Therefore, the formal limit as
$\epsilon\rightarrow 0$ of \eqref{na1}--\eqref{na4} is the following
incompressible MHD equations (suppose that the limits
${\u}^\epsilon\rightarrow \w$ and ${\H}^\epsilon\rightarrow \t$
exist.)
\begin{align}
&\partial_t
\w+\w\cdot\nabla \w+\nabla\pi+\frac{1}{2}\nabla(|\t|^2)-\t\cdot\nabla\t=\mu\Delta\w,\label{im1}\\
&\partial_t\t+\w\cdot\nabla\t-\t\cdot\nabla\w=\nu\Delta\t,\label{im2}\\
&\dv \w=0,\;\;\;\;\dv\t=0.\label{im3}
\end{align}

In this paper we shall establish the above limit rigorously.
Moreover, we shall show that for sufficiently small Mach number, the
compressible flows admit a smooth solution on the time interval where the smooth solution
of the incompressible MHD equations exists. In addition,
we shall also  study the low Mach number limit of the ideal compressible
MHD equations (namely, $\mu=\lambda=\nu=\kappa=0$ in
\eqref{naa}--\eqref{nad}) for which small pressure and entropy
variations are assumed. The convergence rates are obtained in both
cases.

We should point out here that it still remains to be an open problem
to prove rigorously the low Mach number limit of the ideal or full
non-isentropic MHD equations with large temperature variations in
the framework of classical solutions, even in the whole space case,
although the corresponding problems for the non-isentropic Euler and
the full Navier-Stokes equations were solved in the whole space
\cite{MS01,A06} or the bounded domain in \cite{S86}. The reason is
that the presence of the magnetic field and its interaction with
hydrodynamic motion in the MHD flow of large oscillation cause
serious difficulties. We can not apply directly the techniques
developed in \cite{S86,MS01,A06} for the Euler and Navier-Stokes
equations to obtain the uniform estimates for the solutions to the
ideal or full non-isentropic MHD equations. In the present paper,
however, we shall employ an alternative approach, which is based on
the energy estimates for symmetrizable quasilinear
hyperbolic-parabolic systems and the convergence-stability lemma for
singular limit problems \cite{Y1,BY}, to deal with the ideal or full
non-isentropic MHD equations. There are two advantages of this
approach: The first one is that we can rigorously prove the
incompressible limit in the time interval where the limiting system
admits a smooth solution.  The second one is that the estimates we
obtained do not depend on the viscosity and thermodynamic
coefficients, compared with the results in \cite{HL02} where
all-time existence of smooth solutions to the full Navier-Stokes
equations was discussed and  the estimates depend on the parameters
intimately.

For large entropy variation and general initial data, the authors
have rigorously proved the low Mach number limit of the non-isentropic MHD
equations  with zero magnetic diffusivity in \cite{JJL3} by adapting and modifying
the approach developed in \cite{MS01}. We mention that
the coupled singular limit problem for the full MHD equations in the
framework of the so-called variational solutions were studied recently in
\cite{K10,NRT}.

 Before ending the introduction, we give the notations used throughout the current
 paper. We use the letter $C$ to denote various positive constants
independent of $\epsilon$.  For convenience, we denote by $H^l
\equiv H^l (\Omega)$ ($l\in\mathbb{R}$) the standard Sobolev spaces and
 write $\|\cdot\|_l$ for the standard norm of $H^l$ and
$\|\cdot\|$ for $\|\cdot\|_0$.

This paper is organized as follows.
In Section 2 we state our main results. The proof for the full MHD equations
and the ideal MHD equations is presented  in  Section 3 and  Section 4, respectively.
Finally, an appendix is given to derive briefly the dimensionless form of the
full compressible MHD equations.

\section{Main results}
We first recall the local existence of strong solutions to the
incompressible MHD equations \eqref{im1}--\eqref{im3} in the domain
$\Omega$. The proof can be found in \cite{DL,ST}. Recall here that
$\Omega =\mathbb{R}^3 $ or $\Omega =\mathbb{T}^3 $.

\begin{prop}[\!\cite{DL,ST}]\label{imhd}
Let $s>3/2+2$. Assume that the initial data $(\w ,\t )|_{t=0}$ $=(\w_0,\t_0)$ satisfy
$\w_0\in {H}^s, \t_0\in {H}^s$, and $\dv\,\w_0=0$, $\dv\t_0=0$.
Then, there exist a
 $\hat{T}^*\in (0,\infty]$ and a unique solution
 $(\w,\t)\in L^{\infty}(0,\hat{T}^*;{H}^s)$ to the
incompressible MHD equations \eqref{im1}--\eqref{im3}, and  for
any $0<T<\hat{T}^*$,
\begin{equation}   
\sup_{0\le t\le T}\!\big\{||(\w,\t)(t)||_{H^s}
+||(\partial_t\w,\partial_t\t)(t)||_{H^{s-2}}+ ||\nabla
\pi(t)||_{H^{s-2}}\big\} \le C. \nonumber
\end{equation}
\end{prop}

Denoting
$U^\epsilon=(q^\epsilon,\u^\epsilon,\H^\epsilon,\phi^\epsilon)^\top$,
we rewrite the system \eqref{na1}--\eqref{na4} in the vector form
\begin{align}
 A_0(U^\epsilon)\partial_tU^\epsilon+\sum_{j=1}^3A_j(U^\epsilon)\partial_jU^\epsilon=Q(U^\epsilon),\label{sy1}
\end{align}
where
$$
Q(U^\epsilon)  =\big(0,
{F}(\u^\epsilon), \nu\Delta \H^\epsilon, \kappa\Delta
\phi^\epsilon+ \epsilon(L(\u^\epsilon)+G(\H^\epsilon))\big)^\top,
$$
 with
\begin{eqnarray*}
 {F}(\u^\epsilon)&=& 2\mu\dv(\mathbb{D}(\u^\epsilon))+\lambda\nabla(\mbox{tr}\mathbb{D}(\u^\epsilon)),\\
L(\u^\epsilon) &=& 2\mu|\mathbb{D}(\u^\epsilon)|^2+\lambda(\mbox{tr}\mathbb{D}(\u^\epsilon))^2,\\
G(\H^\epsilon) &=& \nu|\nabla\times \H^\epsilon|^2,
\end{eqnarray*}
and the matrices $A_j(U^\epsilon)$ ($0\leq j\leq 3$) are given by
\begin{align*}
& \qquad \qquad \ \ \quad A_0(U^\epsilon)=\mbox{diag}(1,1+\epsilon
q^\epsilon,1+\epsilon q^\epsilon,1+\epsilon
q^\epsilon,1,1,1,1+\epsilon q^\epsilon),\\[1mm]
&\qquad \qquad \ \ \quad  A_1(U^\epsilon)=\\
&{\small \left(
\begin{array}{cccccccc}
u_1^\epsilon&\frac{1+\epsilon q^\epsilon}{\epsilon}&0&0&0&0&0&0\\
\frac{1+\epsilon
\phi^\epsilon}{\epsilon}&u_1^\epsilon(1+\epsilon q^\epsilon)
&0&0&0&H_2^\epsilon&H_3^\epsilon&\frac{1+\epsilon q^\epsilon}{\epsilon}\\
0&0&u_1^\epsilon(1+\epsilon q^\epsilon)&0&0&-H^\epsilon_1&0&0\\
0&0&0&u_1^\epsilon(1+\epsilon q^\epsilon)&0&0&-H_1^\epsilon&0\\
0&0&0&0&u_1^\epsilon&0&0&0\\
0&H_2^\epsilon&-H_1^\epsilon&0&0&u_1^\epsilon&0&0\\
0&H_3^\epsilon&0&-H_1^\epsilon&0&0&u_1^\epsilon&0\\
0&\frac{(\gamma-1)(1+\epsilon q^\epsilon)(1+\epsilon
\phi^\epsilon)}{\epsilon}&0&0&0&0&0&(1+\epsilon
q^\epsilon)u_1^\epsilon
\end{array}
\right),}  \\[1mm]
&\qquad \qquad \ \ \quad A_2(U^\epsilon)=\\
&{\small \left(
\begin{array}{cccccccc}
u_2^\epsilon&0&\frac{1+\epsilon q^\epsilon}{\epsilon}&0&0&0&0&0\\
0&u_2^\epsilon(1+\epsilon q^\epsilon)&0&0&-H_2^\epsilon&0&0&0\\
\frac{1+\epsilon
\phi^\epsilon}{\epsilon}&0&u_2^\epsilon(1+\epsilon q^\epsilon)
&0&H_1^\epsilon&0&H_3^\epsilon&\frac{1+\epsilon q^\epsilon}{\epsilon}\\
0&0&0&u_2^\epsilon(1+\epsilon q^\epsilon)&0&0&-H_2^\epsilon&0\\
0&-H_2^\epsilon&H_1^\epsilon&0&u_2^\epsilon&0&0&0\\
0&0&0&0&0&u_2^\epsilon&0&0\\
0&0&H_3^\epsilon&-H_2^\epsilon&0&0&u_2&0\\
0&0&\frac{(\gamma-1)(1+\epsilon q^\epsilon)(1+\epsilon
\phi^\epsilon)}{\epsilon}&0&0&0&0&(1+\epsilon
q^\epsilon)u_2^\epsilon
\end{array}
\right),}\\[1mm]
&\qquad \qquad \ \ \quad  A_3(U^\epsilon)=\\
&{\small \left(
\begin{array}{cccccccc}
u_3^\epsilon&0&0&\frac{1+\epsilon q^\epsilon}{\epsilon}&0&0&0&0\\
0&u_3^\epsilon(1+\epsilon q^\epsilon)&0&0&-H_3^\epsilon&0&0&0\\
0&0&u_3^\epsilon(1+\epsilon q^\epsilon)&0&0&-H_3^\epsilon&0&0\\
\frac{1+\epsilon
\phi^\epsilon}{\epsilon}&0&0&u_3^\epsilon(1+\epsilon q^\epsilon)
&H_1^\epsilon&H_2^\epsilon&0&\frac{1+\epsilon q^\epsilon}{\epsilon}\\
0&-H_3^\epsilon&0&H_1^\epsilon&u_3^\epsilon&0&0&0\\
0&0&-H_3^\epsilon&H_2^\epsilon&0&u_3^\epsilon&0&0\\0&0&0&0&0&0&u_3^\epsilon&0\\
0&0&0&\frac{(\gamma-1)(1+\epsilon q^\epsilon)(1+\epsilon
\phi^\epsilon)}{\epsilon}&0&0&0&(1+\epsilon q^\epsilon)u_3^\epsilon
\end{array}
\right).}
\end{align*}

 It is easy to see that the matrices
$A_j(U^\epsilon)$ ($0\leq j\leq 3$) can be  symmetrized by choosing
$$\hat{A}_0(U^\epsilon)=\mbox{diag}\big((1+\epsilon
\phi^\epsilon)(1+\epsilon q^\epsilon)^{-1},1,1,1,1,1,1,[(\gamma-1)(1+\epsilon
\phi^\epsilon)]^{-1}\big).$$
Moreover, for $U^\epsilon\in \bar{G}_1\subset\subset G$ with $G$ being the state
space for the system \eqref{sy1}, $\hat{A}_0(U^\epsilon)$ is a positive definite
symmetric matrix for sufficiently small $\epsilon$.

Assume that the initial data
${U}^\epsilon(x,0)={U}_0^\epsilon(x)=({q^\epsilon_0}(x),{\u^\epsilon_0}(x),
{\H^\epsilon_0}(x),{\phi^\epsilon_0}(x))^\top\in
H^{s}$   and ${U}_0^\epsilon(x)\in G_0,$ $
\bar{G}_0\subset\subset G$. The main theorem of the present paper is the following.
\begin{thm}\label{thm1}
Let $s>{3}/{2}+2$. Suppose that the initial data ${U}_0^\epsilon(x)$
satisfy
$$\left\|{U}_0^\epsilon(x)-\left(0,  \w_0(x),\t_0(x),0\right)^\top\right\|_s=O(\epsilon).$$
Let $(\w,\t,\pi)$ be a smooth solution to \eqref{im1}--\eqref{im3}
obtained in Proposition \ref{imhd}. If $(\w, \pi)\in C([0,T^\ast],
H^{s+2})\cap C^1([0,T^\ast],  H^{s})$ with $T^\ast>0$ finite, then
there exists a constant $\epsilon_0>0$ such that, for all $\epsilon\leq
 \epsilon_0$, the system \eqref{sy1} with
initial data ${U}^\epsilon_0(x)$ has a unique smooth solution
${U}^\epsilon(x,t)\in C([0, T^\ast],H^s)$. Moreover,
 there exists a positive constant $K>0$, independent of $\epsilon$, such that, for all
 $\epsilon\leq  \epsilon_0$,
 \begin{align}\label{er1}
\sup_{t\in [0,T^\ast]}\left\|U^\epsilon(\cdot, t)-\left(\frac{\epsilon}{2}\pi,
\w,\t,\frac{\epsilon}{2}\pi\right)^\top\right\|_s\leq K\epsilon.
\end{align}
\end{thm}

\begin{rem}
From Theorem \ref{thm1}, we know that for sufficiently small
$\epsilon$ and well-prepared initial data, the full MHD equations
\eqref{naa}--\eqref{nad} admits a unique smooth solution on the same
time interval where the smooth solution of the incompressible MHD equations exists.
Moreover, the solution can be approximated as shown in \eqref{er1}.
\end{rem}

\begin{rem}
We remark that the constant $K$ in \eqref{er1} is also independent
of the coefficients $\mu, \nu$ and $\kappa$. This is quite different from the results
by Hagstorm and Loranz in \cite{HL02}, where the estimates do
depend on $\mu$  intimately.
\end{rem}

Our approach is still valid for the ideal compressible MHD equations.
However, we will give a particular analysis for the ideal model with
more general pressure by using the entropy form of the energy
equation rather than the thermal energy equation in \eqref{nagg}.

The ideal compressible  MHD equations can be written as
\begin{align}
&\partial_t\rho +\dv(\rho\u)=0, \label{naa2} \\
&\partial_t(\rho\u)+\dv\left(\rho\u\otimes\u\right)+{\nabla
p}
  =\frac{1}{4\pi}(\nabla \times \H)\times \H, \label{nab2} \\
&\partial_t\H-\nabla\times(\u\times\H)=0,\quad \dv\H=0,\label{nac2}\\
&\partial_t\E+\dv\left(\u(\E'+p)\right)
=\frac{1}{4\pi}\dv\big((\u\times\H)\times\H\big).
  \label{nad2}
\end{align}

With the help of the Gibbs relation
\begin{equation*}
\theta \mathrm{d}S=\mathrm{d}e
+p\,\mathrm{d}\left(\frac{1}{\rho}\right)
\end{equation*}
and the  identity \eqref{nae},
the energy balance equation \eqref{nad2}  is replaced by
\begin{equation}\label{naf2}
\partial_t(\rho S)+\dv(\rho  S\u)=0,
\end{equation}
where $S$ denotes the entropy. We   reconsider the equation of state
as a function of $S$ and $p$,   i.e.  $\rho=R(S,p)$ for some
positive smooth function $R$ defined for all $S$ and $p>0$, and
satisfying $\frac{\partial R}{\partial p}>0$. Then, by
  utilizing   \eqref{naa2}, \eqref{idena} and \eqref{idenb},
together with the constraint $\dv {\H}=0$,  the system
\eqref{naa2}--\eqref{nac2}  and \eqref{naf2}  can be
written in the dimensionless form as follows
(after applying the arguments similar to those in the Appendix):
\begin{align}
    & A(S^\epsilon,p^\epsilon)(\partial_t p^\epsilon+\u^\epsilon\cdot \nabla p^\epsilon)+\dv \u^\epsilon=0,\label{nag}\\
& R(S^\epsilon,p^\epsilon)(\partial_t \u^\epsilon+\u^\epsilon\cdot \nabla \u^\epsilon)+\frac{\nabla
p^\epsilon}{\epsilon^2} -\H^\epsilon\cdot\nabla\H^\epsilon+\frac{1}{2}\nabla(|\H^\epsilon|^2)= 0, \label{nah}\\
&  \partial_t {\H^\epsilon}+\u^\epsilon\cdot\nabla\H^\epsilon+\dv \u^\epsilon\H^\epsilon-\H^\epsilon\cdot\nabla
\u^\epsilon=0, \quad \dv \H^\epsilon=0, \label{nai}\\
 &\partial_tS^\epsilon+\u^\epsilon\cdot \nabla S^\epsilon=0,\label{naj}
\end{align}
where $A(S^\epsilon,p^\epsilon)=\frac{1}{R(S^\epsilon,p^\epsilon)}\frac{\partial R(S^\epsilon,p^\epsilon)}{\partial p^\epsilon}$.

To study the low Mach number limit of the above system, we
introduce the transformation
\begin{align}
p^\epsilon (x,t)=\underline{p} e^{\epsilon q^\epsilon(x, t)}, \;\;
S^\epsilon(x, t)=\underline{S}+\epsilon \Theta^\epsilon(x, t),\label{se1}
\end{align}
where $\underline{p}$ and $\underline{S}$ are positive constants, to obtain that
\begin{align}
    & a(\underline{S}+\epsilon\Theta^\epsilon,\epsilon q^\epsilon)(\partial_t q^\epsilon+\u^\epsilon\cdot \nabla q^\epsilon)
    +\frac{1}{\epsilon}\dv \u^\epsilon=0,\label{nak1}\\
& r(\underline{S}+\epsilon\Theta^\epsilon,\epsilon q^\epsilon)(\partial_t
\u^\epsilon+\u^\epsilon\cdot \nabla \u^\epsilon)+\frac{1}{\epsilon}\nabla q^\epsilon
 -\H^\epsilon\cdot\nabla\H^\epsilon+\frac{1}{2}\nabla(|\H^\epsilon|^2)= 0,  \label{nal1}\\
&  \partial_t {\H}^\epsilon +\u^\epsilon\cdot\nabla\H^\epsilon+\dv \u^\epsilon\H^\epsilon-\H^\epsilon\cdot\nabla
\u^\epsilon=0, \quad \dv \H^\epsilon=0, \label{nam1}\\
 &\partial_t\Theta^\epsilon+\u^\epsilon\cdot \nabla\Theta^\epsilon=0,   \label{nan1}
\end{align}
where
\begin{gather*}
 a(S^\epsilon,\epsilon q^\epsilon)= A(S^\epsilon, \underline{p}e^{\epsilon
q^\epsilon})\underline{p}e^{\epsilon q^\epsilon}
=\frac{\underline{p}e^{\epsilon
q^\epsilon}}{R(S^\epsilon,\underline{p}e^{\epsilon
 q^\epsilon})}\cdot
 \frac{\partial R(S^\epsilon,s)}{\partial s}\Big|_{s=\underline{p}e^{\epsilon q^\epsilon}},\\ 
  r^\epsilon(S^\epsilon,\epsilon q^\epsilon) = \frac{R(S^\epsilon,\underline{p}e^{\epsilon
 q^\epsilon})}{\underline{p}e^{\epsilon q^\epsilon}}.   
 \end{gather*}

Making use of the fact that $\cu \nabla =0$ and letting
  $\epsilon \rightarrow 0$ in \eqref{nak1} and \eqref{nal1}, we formally deduce that
$\dv \v=0$ and
\begin{align*}
\cu\big( r(\underline{S},0)(\partial_t \v+\v\cdot \nabla \v)
 -(\nabla \times{\d}) \times{\d} \big)=0,
\end{align*}
where we have supposed that the limits  ${\u}^\epsilon\rightarrow \v$ and
${\H}^\epsilon\rightarrow \d$ exist.  Thus,  we can
expect that the limiting system of \eqref{nak1}--\eqref{nan1} takes
the form
\begin{align}
&   r(\underline{S},0)(\partial_t \v+\v\cdot \nabla \v)
  -(\nabla \times{\d}) \times {\d} +\nabla \Pi =0, \label{nao1} \\
&  \partial_t {\d}  +  {\v}  \cdot \nabla {\d}
   -  {\d} \cdot \nabla {\v}  =0, \label{nap1} \\
& \dv \v=0,  \quad \dv {\d} =0  \label{nar1}
\end{align}
for some function $\Pi$.

In order to state our result, we first recall the local existence of
strong solutions to the ideal incompressible MHD equations
\eqref{nao1}--\eqref{nar1} in the domain $\Omega$. The proof can be
found in \cite{DL,ST}.

\begin{prop}[\!\cite{DL,ST}]\label{ivmhd}
Let $s>3/2+1$. Assume that the initial data $(\v ,\d )|_{t=0}$ $=(\v_0,\d_0)$ satisfy
$\v_0\in {H}^s, \d_0\in {H}^s$, and $\dv\,\v_0=0$, $\dv\d_0=0$.
Then, there exist a
 $\tilde{T}^*\in (0,\infty]$ and a unique smooth solution
 $(\v,\d)\in L^{\infty}(0,\tilde{T}^*; {H}^s)$ to the
incompressible MHD equations \eqref{im1}--\eqref{im3}, and  for
any $0<T<\tilde{T}^*$,
\begin{equation}
\sup_{0\le t\le T}\!\big\{||(\v,\d)(t)||_{H^s}
+||(\partial_t\v,\partial_t\d)(t)||_{H^{s-1}}+ ||\nabla
\Pi(t)||_{H^{s-1}}\big\} \le C. \nonumber
\end{equation}
\end{prop}

In the vector form, we arrive at,
for $V^\epsilon=(q^\epsilon,\u^\epsilon,\H^\epsilon,\Theta^\epsilon)^\top$, that
\begin{align}
A_0(\epsilon\Theta^\epsilon,\epsilon q^\epsilon)\partial_t
V^\epsilon+\sum_{j=1}^3\Big\{u^\epsilon_jA_0(\epsilon\Theta^\epsilon,\epsilon
q^\epsilon)+\epsilon^{-1}C_j+B_j(\H^\epsilon)\Big\}\partial_jV^\epsilon=0,\label{syin}
\end{align}
where
$$A_0(\epsilon\Theta^\epsilon,\epsilon q^\epsilon)=\mbox{diag}(a(S^\epsilon,\epsilon q^\epsilon),
r(S^\epsilon,\epsilon q^\epsilon),r(S^\epsilon,\epsilon q^\epsilon),r(S^\epsilon,\epsilon q^\epsilon),1,1,1,1),
$$
and $C_j$ is symmetric constant  matrix,
 and $B_j(\H^\epsilon)$ is a symmetric matrix of $\H^\epsilon$.

 Assume that the initial data for the equations \eqref{syin} satisfy
$${V}^\epsilon_0(x)=(\tilde{q}_0^\epsilon(x),\tilde{\u}_0^\epsilon(x),\tilde{\H}_0^\epsilon(x),
\tilde{\Theta}_0^\epsilon(x))^\top\in
H^{s},\;\mbox{ and }\;V_0^\epsilon(x)\in G_0,\;\;
\bar{G}_0\subset\subset G$$ with $G$ being state space for
\eqref{syin}.
Thus, our result on the ideal compressible MHD equations reads as
\begin{thm}\label{thm2}
Let $s>{3}/{2}+1$. Suppose that the initial data ${V}_0^\epsilon(x)$
satisfy
$$\left\|{V}_0^\epsilon(x)-\left(0,  {\v}_0(x), \d_0(x),0\right)^\top\right\|_s=O(\epsilon).$$
Let $(\v,\d,\Pi)$ be a smooth solution to \eqref{nao1}--\eqref{nar1}  obtained in
Proposition \ref{ivmhd}. If $(\v,\Pi)\in C([0,\bar T_\ast],
H^{s+1})\cap C^1([0,\bar T_\ast]$, $H^s)$ with $\bar T_\ast>0$ finite,
then there exists  a constant $\epsilon_1>0$ such that, for all $\epsilon\leq
 \epsilon_1$, the system \eqref{syin} with
initial data ${V}_0^\epsilon(x)$ has a unique solution
$V^\epsilon(x, t)\in C([0,\bar
T_\ast],H^s)$. Moreover,
 there exists a positive constant $K_1>0$ such that, for all $\epsilon\leq
 \epsilon_1$,
 \begin{align}\label{er2}
\sup_{t\in [0,\bar T_\ast]}\left\|V^\epsilon(\cdot, t)-(\epsilon\Pi,
\v,\d,\epsilon\Pi)^\top\right\|_s\leq K_1\epsilon.
\end{align}
\end{thm}

\section{Proof of Theorem \ref{thm1}}

This section is devoted to proving Theorem \ref{thm1}. First, following the
proof of the local existence theory for the initial value problem of
symmetrizable hyperbolic-parabolic systems by Volpert and Hudjaev in \cite{VH}, we
obtain that there exists a time interval $[0,T]$ with $T>0$, so that
the system \eqref{sy1} with initial data ${U}_0^\epsilon(x)$ has a unique classical
solution $U^\epsilon(x,t)\in C([0,T],H^s)$
and $U^\epsilon(x,t)\in G_2$ with $\bar{G}_2\subset\subset G$. We remark that
the crucial step in the proof  of local existence result is to
prove the uniform boundedness of the solutions. See also \cite{JLL} for some relative results.

Now, define
\begin{align*}
T_\epsilon=\sup\{T>0: U^\epsilon(x,t)\in C([0,T],H^s), U^\epsilon(x,t)\in G_2,
\forall \, (x,t)\in \Omega\times [0,T]\}.
\end{align*}
Note that $T_\epsilon$ depends on $\epsilon$ and may tend to
zero as $\epsilon$ goes to $0$.

To show that $\underline{\lim}_{\epsilon\rightarrow
0}T_\epsilon
>0$,  we shall make use of the
convergence-stability lemma which was established in  \cite{Y1,BY} for hyperbolic systems of balance laws. It is also implied in \cite{Y1} that a convergence-stability lemma can be formulated as a part of (local) existence theories for any evolution equations. For the hyperbolic-parabolic system \eqref{sy1}, we have the following  convergence-stability lemma.
\begin{lem}\label{lem1}
Let   $s>{3}/{2}+2$. Suppose that ${U}^\epsilon_0(x)\in G_0, \bar{G}_0\subset\subset G,$ and
${U}^\epsilon_0(x)\in H^{s}$, and the following
convergence assumption (A) holds.

(A) There exists $T_\star>0$ and $U_\epsilon\in
L^\infty(0,T_\star ;H^s)$ for each $\epsilon$, satisfying
\begin{align}
\bigcup_{x,t,\epsilon}\{U_\epsilon(x,t)\}\subset\subset G,\nonumber
\end{align}
such that for $t\in [0, \min\{T_\star,T_\epsilon\})$,
$$
\sup_{x,t}|U^\epsilon(x,t)-U_{\epsilon}(x,t)|=o(1),\;\;
\sup_{t}\|U^\epsilon(x,t)-U_{\epsilon}(x,t)\|_s=O(1),\quad\mbox{as }\epsilon\to 0.
$$
Then, there exist an $\bar \epsilon>0$ such that, for all $\epsilon\in (0,\bar \epsilon]$, it holds that
\begin{align}
T_\epsilon>T_\star.\nonumber
\end{align}

\end{lem}

To apply Lemma \ref{lem1}, we construct the approximation
$U_\epsilon=(q_\epsilon,\v_\epsilon,\t_\epsilon,\phi_\epsilon)^\top$
with $q_\epsilon=\epsilon\pi/2, \v_\epsilon=\w,
\t_\epsilon=\t,$ and $\phi_\epsilon={\epsilon}\pi/2$,
where $(\w,\t,\pi)$ is the classical solution to the system \eqref{im1}--\eqref{im3}
obtained in Proposition \ref{imhd}. It is easy to verify that
$U_\epsilon$
satisfies
\begin{align}
&\partial_tq_\epsilon+\v_\epsilon\cdot\nabla
q_\epsilon+\frac{1}{\epsilon}(1+\epsilon
q_\epsilon)\dv\v_\epsilon=
\frac{\epsilon}{2}(\pi_t+\w\cdot\nabla\pi),\label{ap1}\\
&(1+\epsilon
q_\epsilon)(\partial_t\v_\epsilon+\v_\epsilon\cdot\nabla
\v_\epsilon)+\frac{1}{\epsilon}\big[(1+\epsilon
q_\epsilon)\nabla \phi_\epsilon+(1+\epsilon
\phi_\epsilon)\nabla q_\epsilon\big]\nonumber\\
&\qquad\qquad -\t_\epsilon\cdot\nabla\t_\epsilon +\frac{1}{2}\nabla(|\t_\epsilon|^2)=\mu\Delta
\v_\epsilon+\frac{\epsilon^2}{2}\pi(\w_t+\w\cdot\nabla\w+\nabla\pi),\label{ap2}\\[1mm]
&\partial_t\t_\epsilon+\v_\epsilon\cdot\nabla\t_\epsilon+
\dv \v_\epsilon\t_\epsilon-\t_\epsilon\cdot\nabla
\v_\epsilon=\nu\Delta\t_\epsilon,\quad
\dv\t_\epsilon=0,\label{ap3}\\
&(1+\epsilon q_\epsilon)(\partial_t
\phi_\epsilon+\v_\epsilon\cdot\nabla
\phi_\epsilon)+\frac{\gamma-1}{\epsilon}(1+\epsilon
q_\epsilon)(1+\epsilon \phi_\epsilon)\dv
\v_\epsilon\nonumber\\
& \qquad \qquad \qquad\qquad\qquad\qquad\qquad
 =\left(\frac{\epsilon}{2}+\frac{\epsilon^3}{4}\pi\right)(\pi_t+\w\cdot\nabla\pi). \label{ap4}
\end{align}

We rewrite the system \eqref{ap1}--\eqref{ap4} in the following vector form
\begin{align}
A_0(U_\epsilon)\partial_t
U_\epsilon+\sum_{j=1}^3A_j(U_\epsilon)\partial_jU_\epsilon=S(U_\epsilon)+R,\label{sy2}
\end{align}
with $S(U_\epsilon)=(0,\mu\Delta \v_\epsilon,\nu\Delta\t_\epsilon,0)^\top$ and
\begin{equation*}
R=\left(
\begin{array}{c}
\frac{\epsilon}{2}(\pi_t+\w\cdot\nabla\pi)\\
\frac{\epsilon^2}{2}\pi(\w_t+\w\cdot\nabla\w+\nabla\pi)\\
\big(\frac{\epsilon}{2}+\frac{\epsilon^3}{4}\pi\big)(\pi_t+\w\cdot\nabla\pi)\\
0 \end{array} \right).
\end{equation*}
Due to the regularity assumptions on $(\w,\pi)$ in Theorem
\ref{thm1}, we have
\begin{align}
\max_{t\in [0,T^\ast]}\|R(t)\|_s\leq C\epsilon. \nonumber
\end{align}

To prove Theorem \ref{thm1},  it suffices to prove the error
estimate in \eqref{er1} for $t\in [0,\min\{T^\ast, T_\epsilon\})$
thanks to Lemma \ref{lem1}. To this end, introducing
$$E=U^\epsilon-U_\epsilon\quad\mbox{ and }\quad
\mathcal{A}_j(U)=A_0^{-1}(U)A_j(U),$$
and using \eqref{sy1} and \eqref{sy2}, we see that
\begin{align}
E_t+\sum_{j=1}^3\mathcal{A}_j(U^\epsilon)E_{x_j}=&(\mathcal{A}_j(U_\epsilon)-\mathcal{A}_j(U^\epsilon))
U_{\epsilon x_j}+A_0^{-1}(U^\epsilon)Q(U^\epsilon)\nonumber\\
& -A_0^{-1}(U_\epsilon)(S(U_\epsilon)+R).\label{E1}
\end{align}

For any multi-index $\alpha$ satisfying $|\alpha|\leq s$, we take
the operator $D^\alpha$ to \eqref{E1} to obtain
\begin{align}
\partial_tD^\alpha E+\sum_{j=1}^3\mathcal{A}_j(U^\epsilon)\partial_{x_j}D^\alpha
E=P_1^\alpha+P_2^\alpha+Q^\alpha+R^\alpha\label{EA1}
\end{align}
with
\begin{align}
&P_1^\alpha=\sum_{j=1}^3\{\mathcal{A}_j(U^\epsilon)\partial_{x_j}D^\alpha
E-D^\alpha(\mathcal{A}_j(U^\epsilon)\partial_{x_j}
E)\},\nonumber\\
&P_2^\alpha=\sum_{j=1}^3D^\alpha\{(\mathcal{A}_j(U_\epsilon)-\mathcal{A}_j(U^\epsilon))
U_{\epsilon
x_j}\},\nonumber\\
&Q^\alpha=D^\alpha\{A_0^{-1}(U^\epsilon)Q(U^\epsilon)-A_0^{-1}(U_\epsilon)S(U_\epsilon)\},
\nonumber\\[1mm]
& R^\alpha=D^\alpha\{A_0^{-1}(U_\epsilon)R\}.\nonumber
\end{align}

Define
$$ \tilde{A}_0(U^\epsilon)=\mbox{diag}\Big(\frac{1+\epsilon \phi^\epsilon}
{(1+\epsilon q^\epsilon)^2},1,1,1,\frac{1}{1+\epsilon q^\epsilon},
\frac{1}{1+\epsilon q^\epsilon}, \frac{1}{1+\epsilon q^\epsilon},
\frac{1}{(\gamma-1)(1+\epsilon \phi^\epsilon)}\Big), $$
and the canonical energy by
\begin{align}
\|E\|_\mathrm{e}^2:=\int \langle\tilde{A}_0(U^\epsilon) E, E\rangle dx.\nonumber
\end{align}
Note that $\tilde{A}_0(U^\epsilon)$ is a positive definite symmetric matrix
and $\tilde{A}_0(U^\epsilon)\mathcal{A}_j(U^\epsilon)$ is symmetric. Now, if we multiply
\eqref{EA1} with $\tilde{A}_0(U^\epsilon)$ and take the inner product
between the resulting system and $D^\alpha E$, we arrive at
\begin{align}
\frac{d}{dt}\|D^\alpha E\|_{\mathrm{{e}}}^2= & 2\int \langle \Gamma D^\alpha E,
D^\alpha E\rangle dx\nonumber\\
& + 2\int (D^\alpha E)^T\tilde{A}_0(U^\epsilon)(P_1^\alpha+P_2^\alpha+Q^\alpha+R^\alpha),
\label{de1}
\end{align}
where
\begin{align}
\Gamma=(\partial_t,\nabla)\cdot\Big(\tilde{A}_0,\tilde{A}_0(U^\epsilon)\mathcal{A}_1(U^\epsilon),
\tilde{A}_0(U^\epsilon)\mathcal{A}_2(U^\epsilon),\tilde{A}_0(U^\epsilon)\mathcal{A}_3(U^\epsilon)\Big).
\nonumber
\end{align}

Next, we estimate various terms on the right-hand side of  \eqref{de1}. Note that our
estimates only need to be done for $t\in
[0,\min\{T^\ast,T_\epsilon\})$, in which both $U^\epsilon$ and
$U_\epsilon$ are regular enough and take values in a convex
compact subset of the state space. Thus, we have
\begin{align}\label{energy1}
C^{-1}\int |D^\alpha E|^2 \leq \|D^\alpha E\|_{\mathrm{e}}^2\leq C\int |D^\alpha E|^2
\end{align}
and
\begin{align}
|(D^\alpha E)^\top\tilde{A}_0(U^\epsilon)(P_1^\alpha+P_2^\alpha+R^\alpha)| \leq
C(|D^\alpha E|^2+|P_1^\alpha|^2+|P_2^\alpha|^2+|R^\alpha|^2).
\nonumber
\end{align}

To estimate $\Gamma$, we write
$\mathcal{A}_j(U^\epsilon)=u^\epsilon_j\mathbf{I}_{8}+\bar{\mathcal{A}}_j(U^\epsilon)$. Notice that
$\bar{\mathcal{A}}_j(U^\epsilon)$ depends only on $q^\epsilon, \phi^\epsilon$ and $ \H^\epsilon$.
Thus using \eqref{na1} and \eqref{na4}, we have
\begin{align}
|\Gamma|=& \left|\frac{\partial}{\partial
t}\tilde{A}_0+u_j^\epsilon\cdot\nabla\tilde{A}_0+\tilde{A}_0\dv
\u^\epsilon+\dv(\tilde{A}_0\bar{\mathcal{A}}_j(U^\epsilon))\right|\nonumber\\
=& \big| \tilde{A}_0\,\dv \u^\epsilon-\tilde{A}_{0\eta_1}^{'}(1+\epsilon q^\epsilon)\dv
\u^\epsilon-\tilde{A}_{0\eta_2}^{'}[(1+\epsilon \phi^\epsilon)\dv
\u^\epsilon+\kappa(1+\epsilon q^\epsilon)^{-1} \Delta \phi^\epsilon\nonumber\\
& +\epsilon^2(L(\u^\epsilon)+G(\H^\epsilon))]
+ \dv(\tilde{A}_0\bar{\mathcal{A}}_j(U^\epsilon))\big|\nonumber\\
\leq & C+C(|\nabla\u^\epsilon|+|\nabla q^\epsilon|+|\nabla
\phi^\epsilon|
+|\nabla \H^\epsilon|+|\Delta \phi^\epsilon|+|\nabla\u^\epsilon|^2+|\nabla \H^\epsilon|^2)\nonumber\\
\leq&  C+C(|\nabla E|+|\nabla E|^2)+C|\Delta (\phi^\epsilon-\phi_\epsilon)|
+C(|\nabla U_\epsilon|+|\nabla U_\epsilon|^2)\nonumber\\
\leq &  C+C(\|E\|_s+\|E\|_s^2),   \nonumber 
\end{align}
where we have used Sobolev's embedding theorem and the fact that
$s>{3}/{2}+2$, and the symbols $\tilde{A}_{0\eta_1}^{'}$ and $\tilde{A}_{0\eta_2}^{'}$ denote
the differentiation of $\tilde{A}_{0}$ with respect to $\rho^\epsilon$ and $\theta^\epsilon$,
respectively.

Since
\begin{align*}
&\mathcal{A}_j(U^\epsilon)\partial_{x_j}D^\alpha
E-D^\alpha(\mathcal{A}_j(U^\epsilon)\partial_{x_j}
E)= -\sum_{0<\beta\leq\alpha} \binom{\alpha}{\beta}
\partial^\beta\mathcal{A}_j(U^\epsilon)\partial^{\alpha-\beta}E_{x_j}\\
& \qquad = -\sum_{0<\beta\leq\alpha} \binom{\alpha}{\beta}
\partial^\beta[u^\epsilon_j\mathbf{I}_{8}+\bar{\mathcal{A}}_j(U^\epsilon)]\partial^{\alpha-\beta}E_{x_j},
\end{align*}
we obtain, with the help of the Moser-type calculus inequalities in Sobolev
spaces, that
\begin{align*}
\|P_1^\alpha\|\leq &  C\left\{(1+\|(\u^\epsilon,
\H^\epsilon)\|_s)\|E_{x_j}\|_{|\alpha|-1}+\|\epsilon^{-1}(\partial^\beta
f(q^\epsilon,\phi^\epsilon)\partial^{\alpha-\beta}E_{x_j})\|\right\}\nonumber\\
& +C\|\partial^\beta[(1+\epsilon
q^\epsilon)^{-1}(\H^\epsilon-\t_\epsilon)+((1+\epsilon q^\epsilon)^{-1}-(1+\epsilon
q_\epsilon)^{-1})\t_\epsilon]\partial^{\alpha-\beta}E_{x_j}\|\nonumber\\
\leq & C(1+\|E\|_{s}+\|(q^\epsilon,
\phi^\epsilon)\|_s^s)\|E_{x_j}\|_{|\alpha|-1}\nonumber\\
\leq & C(1+\|E\|_s^s)\|E_{x_j}\|_{|\alpha|},
\end{align*}
where $f(q^\epsilon,\phi^\epsilon)=(1+\epsilon q^\epsilon)+(\gamma-1)(1+\epsilon \phi^\epsilon)+
(1+\epsilon q^\epsilon)^{-1}(1+\epsilon \phi^\epsilon).$

Similarly, utilizing the boundedness of $\|U_\epsilon\|_{s+1}$, the term
$P_2^\alpha$ can be bounded as follows:
\begin{align*}
\|P_2^\alpha\|&\leq C\|U_{\epsilon
x_j}\|_s\|\mathcal{A}_j(U_\epsilon)-\mathcal{A}_j(U^\epsilon)\|_{|\alpha|}\nonumber\\
&\leq C\|(u_j^\epsilon-v_{\epsilon
j})\mathbf{I}_{8}+\bar{\mathcal{A}}_j(U^\epsilon)-\bar{\mathcal{A}}_j(U_\epsilon)\|_{|\alpha|}\nonumber\\
&\leq
C(1+\|\u^\epsilon-\v_\epsilon\|_{|\alpha|}+\|\H^\epsilon-\t_\epsilon\|_{|\alpha|})
+C\|\epsilon^{-1}(f(q^\epsilon,\phi^\epsilon)-f(q_\epsilon,\phi_\epsilon))\|_{|\alpha|}\nonumber\\
&\leq C(1+\|q_\epsilon+\eta_3(q^\epsilon-q_\epsilon)+\phi_\epsilon+
\eta_4(\phi^\epsilon-\phi_\epsilon)\|_s^s)\|E\|_{|\alpha|}\nonumber\\
&\leq  C(1+\|E\|_{s}^s)\|E\|_{|\alpha|},
\end{align*}
where $0\leq\eta_3,\eta_4\leq 1$ are  constants.

The estimate of $\int (D^\alpha E)^\top\tilde{A}_0(U^\epsilon)Q^\alpha$ is more
complex and delicate.  First, we can rewrite
$\int (D^\alpha E)^\top\tilde{A}_0(U^\epsilon)Q^\alpha$ as
\begin{align*}
& \int (D^\alpha E)^\top\tilde{A}_0(U^\epsilon)Q^\alpha
=   \int D^\alpha(\u^\epsilon-\v_\epsilon)D^\alpha\big[(1+\epsilon
q^\epsilon)^{-1}F(\u^\epsilon)-\mu(1+\epsilon
q_\epsilon)^{-1}\Delta\v_\epsilon\big]\nonumber\\
&\quad + \nu\int D^\alpha(\H^\epsilon-\t_\epsilon)(1+\epsilon
q^\epsilon)^{-1} D^\alpha(\Delta \H^\epsilon- \Delta
\t_\epsilon)\nonumber\\
&\quad + \kappa(\gamma-1)^{-1}\int D^\alpha(\phi^\epsilon-\phi_\epsilon) (1+\epsilon
\phi^\epsilon)^{-1}D^\alpha\{(1+\epsilon q)^{-1}\Delta \phi^\epsilon-(1+\epsilon
q_\epsilon)^{-1} \Delta
\phi_\epsilon\} \\
& \quad +\epsilon (\gamma-1)^{-1} \int D^\alpha(\phi^\epsilon-\phi_\epsilon)(1+\epsilon
\phi^\epsilon)^{-1}D^\alpha\{(1+\epsilon
q^\epsilon)^{-1} (L(\u^\epsilon)+G(\H^\epsilon))\}\nonumber\\
& = \mathcal{Q}_{u}+\mathcal{Q}_{H}+\mathcal{Q}_{\phi_1} +\mathcal{Q}_{\phi_2}.
\end{align*}
By integration by parts, the Cauchy and Moser-type inequalities, and
Sobolev's embedding theorem, we find that $\mathcal{Q}_{u}$ can be controlled
as follows:
\begin{align*}
\mathcal{Q}_u  = &\int D^\alpha(\u^\epsilon-\v_\epsilon)D^\alpha\{(1+\epsilon
q^\epsilon)^{-1} \mu\Delta(\u^\epsilon-\v_\epsilon)+(\mu+\lambda)\nabla\dv(\u^\epsilon-\v_\epsilon)\}\nonumber\\
&+\mu\int D^\alpha(\u^\epsilon-\v_\epsilon)D^\alpha\{[(1+\epsilon
q^\epsilon)^{-1}-(1+\epsilon q_\epsilon)^{-1}]\Delta
\v_\epsilon\}\nonumber\\
\leq & -\int\frac{\mu}{1+\epsilon
q^\epsilon}|D^\alpha\nabla(\u^\epsilon-\v_\epsilon)|^2-\int\frac{\mu+\lambda}{1+\epsilon
q^\epsilon}|D^\alpha\dv(\u^\epsilon-\v_\epsilon)|^2\nonumber\\
& +\int
D^\alpha(\u^\epsilon-\v_\epsilon)\sum_{0<\beta\leq\alpha}D^\beta[(1+\epsilon
q^\epsilon)^{-1}]D^{\alpha-\beta}\big\{\mu\Delta(\u^\epsilon-\v_\epsilon)\nonumber\\
& +(\mu+\lambda)\nabla\dv(\u^\epsilon-\v_\epsilon)\big\}+C\|E\|_{|\alpha|}^2
+C\|E\|_{s}^4+C\epsilon\|D^{\alpha}\nabla(\u^\epsilon-\v_\epsilon)\|^2\nonumber\\
\leq & -C\int\mu|D^\alpha\nabla(\u^\epsilon-\v_\epsilon)|^2-C\int(\mu+\lambda)
|D^\alpha\dv(\u^\epsilon-\v_\epsilon)|^2+C\|E\|_{|\alpha|}^2\nonumber\\
&+C\epsilon\|D^{\alpha}\nabla(\u^\epsilon-\v_\epsilon)\|^2
+C\|E\|_{s}^4 +C\|E\|_s^2\|E\|_{|\alpha|}^2+
\int D^\alpha(\u^\epsilon-\v_\epsilon) \cdot \nonumber\\
&\ \ \ \ \sum_{1<\beta\leq\alpha} D^\beta[(1+\epsilon
q^\epsilon)^{-1}]D^{\alpha-\beta}\{\mu\Delta(\u^\epsilon-\v_\epsilon)
+(\mu+\lambda)\nabla\dv(\u^\epsilon-\v_\epsilon)\}\nonumber\\
\leq & -C\int\mu|D^\alpha\nabla(\u^\epsilon-\v_\epsilon)|^2
-C\int(\mu+\lambda)|D^\alpha\dv(\u^\epsilon-\v_\epsilon)|^2\nonumber\\
& +C\epsilon\|D^{\alpha}\nabla(\u^\epsilon-\v_\epsilon)\|^2+C\|E\|_s^4+C\|E\|_{|\alpha|}^2.
\end{align*}

Similarly, the terms $\mathcal{Q}_H$, $\mathcal{Q}_{\phi_1}$ and
$\mathcal{Q}_{\phi_2}$ can be bounded as follows:
\begin{eqnarray*} &&
\mathcal{Q}_H\leq -C\nu\int|D^\alpha\nabla(\H^\epsilon-\t_\epsilon)|^2+
C\epsilon\|D^{\alpha}\nabla(\H^\epsilon-\t_\epsilon)\|^2+C\|E\|_s^4 +C\|E\|_{|\alpha|}^2,
\\
&& \mathcal{Q}_{\phi_1}\leq -C\kappa\int|D^\alpha\nabla(\phi^\epsilon-\phi_\epsilon)|^2
+ C\epsilon\|D^{\alpha}\nabla(\phi^\epsilon-\phi_\epsilon)\|^2+C\|E\|_s^4
+C\|E\|_{|\alpha|}^2
\end{eqnarray*}
and
$$ \mathcal{Q}_{\phi_2}\leq C\epsilon\|D^{\alpha}
\nabla(\phi^\epsilon-\phi_\epsilon)\|^2+C\|E\|_s^4 +C\|E\|_{|\alpha|}^2.
$$

 Putting all the above estimates into \eqref{de1} and taking
$\epsilon$ small enough, we obtain that
\begin{align}
& \frac{d}{dt}\|D^\alpha E\|_{\mathrm{e}}^2
+\xi\int|D^\alpha\nabla(\u^\epsilon-\v_\epsilon)|^2
+\nu\int|D^\alpha\nabla(\H^\epsilon-\t_\epsilon)|^2\nonumber\\
&\qquad +\kappa\int
|D^\alpha\nabla(\phi^\epsilon-\phi_\epsilon)|^2\leq
C\|R^\alpha\|^2+C(1+\|E\|_s^{2s})\|E\|_{|\alpha|}^2+\|E\|_s^{4},
\label{de2}
\end{align}
where we have used the following estimate
\begin{align}
\mu\int|D^\alpha\nabla(\u^\epsilon-\v_\epsilon)|^2+(\mu+\lambda)\int|D^\alpha\dv(\u^\epsilon-\v_\epsilon)|^2
\geq \xi\int|D^\alpha\nabla(\u^\epsilon-\v_\epsilon)|^2\nonumber
\end{align}
for some positive constant $\xi>0$.

 Using \eqref{energy1}, we
integrate the inequality \eqref{de2} over $(0,t)$ with
$t<\min\{T_\epsilon, T^\ast\}$ to obtain
\begin{align}
\|D^\alpha E(t)\|^2\leq &  C\|D^\alpha E(0)\|^2+C\int_0^t\|R^\alpha(\tau)\|^2
d\tau\nonumber\\
& +C\int_0^t\big\{(1+\|E\|_s^{2s})\|E\|_{|\alpha|}^2+\|E\|_s^{4}\big\}(\tau)d\tau.
\nonumber   \end{align}
Summing up this inequality for all $\alpha$ with $|\alpha|\leq s$, we get
\begin{align}
\| E(t)\|_s^2\leq C\| E(0)\|_s^2+C\int_0^{T^\ast}\|R(\tau)\|_s^2 d\tau
+C\int_0^t\big\{(1+\|E\|_s^{2s})\|E\|_s^2\big\}(\tau)d\tau.\nonumber
\end{align}
With the help of Gronwall's lemma and the fact that
 $$\| E(0)\|_s^2+\int_0^{T^\ast}\|R(t)\|_s^2
 dt=O(\epsilon^2),$$
 we conclude that
\begin{align}
\| E(t)\|_s^2\leq
C\epsilon^2\mbox{exp}\left\{C\int_0^t(1+\|E(\tau)\|_s^{2s})d\tau\right\}\equiv
\Phi(t).\nonumber
\end{align}
It is easy to see that $\Phi(t)$ satisfies
\begin{align}
\Phi'(t)=C(1+\|E(t)\|_s^{2s})\Phi(t)\leq C\Phi(t)
+C\Phi^{s+1}(t).\nonumber
\end{align}
Thus, employing the nonlinear Gronwall-type inequality, we conclude that
there exists a constant $K$, independent of $\epsilon$, such that
\begin{align}
\|E(t)\|_s\leq K \epsilon,\nonumber
\end{align}
for all $t\in [0,\min\{T_\epsilon,T^*\})$, provided
$\Phi(0)=C\epsilon^2<\mbox{exp}(-CT^*)$.  Thus, the proof is completed.

\section{Proof of Theorem \ref{thm2}}

The proof of Theorem \ref{thm2} is essentially similar to that of Theorem
\ref{thm1}, and we only give some explanations here. The local existence
of classical solution to the system \eqref{syin} is given by the
proof of Theorem 2.1 in \cite{M84}. For each fixed $\epsilon$, we assume that the
maximal time interval of existence is $[0, T^\epsilon)$. To prove
Theorem \ref{thm2}, it is crucial to obtain the error estimates in
\eqref{er2}. For this purpose, we construct the approximation
$V_\epsilon=(q_\epsilon,\v_\epsilon,\d_\epsilon,\Theta_\epsilon)^\top$
with $q_\epsilon={\epsilon}\Pi, \v_\epsilon=\v, \d_\epsilon=\d,$ and
$\Theta_\epsilon=\epsilon\Pi$. It is then easy to verify that
$V_\epsilon$ satisfies
\begin{align}
&a(\underline{S}+\epsilon \Theta_\epsilon, \epsilon
q_\epsilon) (\partial_tq_\epsilon+\v_\epsilon\cdot\nabla
q_\epsilon)+\frac{1}{\epsilon} \dv\v_\epsilon\nonumber\\
&\qquad\qquad \qquad\qquad= \epsilon
a(\underline{S}+\epsilon^2 \Pi, \epsilon^2
\Pi)(\Pi_t+\v\cdot\nabla\Pi),\label{iap1}\\
&r(\underline{S}+\epsilon \Theta_\epsilon, \epsilon
q_\epsilon)(\partial_t\v_\epsilon+\v_\epsilon\cdot\nabla
\v_\epsilon)+\frac{1}{\epsilon}\nabla
q_\epsilon-\d_\epsilon\cdot\nabla \d_\epsilon+\frac{1}{2}(|\d_\epsilon|^2)\nonumber\\
&\qquad\qquad\qquad\qquad=[r(\underline{S}+\epsilon \Theta_\epsilon,
\epsilon q_\epsilon)-r(\underline{S},0)]
(\v_t+\v\cdot\nabla\v),\label{iap2}\\[1mm]
&\partial_t\d_\epsilon+\v_\epsilon\cdot\nabla\d_\epsilon+
\dv \v_\epsilon\d_\epsilon-\d_\epsilon\cdot\nabla
\v_\epsilon=0,\quad \dv\d_\epsilon=0,\label{iap3}\\
&\partial_t\Theta_\epsilon+\v_\epsilon\cdot\nabla
\Theta_\epsilon=\epsilon (\Pi_t+\v\cdot\nabla\Pi).\label{iap4}
\end{align}
Thus we can rewrite \eqref{iap1}--\eqref{iap4} in the vector form of
\eqref{syin} with a source term.
  Letting $E=V^\epsilon-V_\epsilon$, we can perform the energy estimates similar
to those in the proof of Theorem \ref{thm1} to show Theorem \ref{thm2}. Here
we omit the details of the proof for conciseness.

\section{Appendix}
We give a dimensionless form of the system \eqref{naa}-\eqref{nac} and \eqref{nagg}
for the ionized fluid obeying the perfect gas relations \eqref{pg} by following
the spirit of \cite{Hu}. Introduce the new dimensionless quantities:
\begin{gather}
 {x}_\star=\frac{{x}}{L_0},\;\;t_\star=\frac{t}{L_0/u_0},\;\;\u_\star=\frac{\u}{u_0},\nonumber\\
 \H_\star=\frac{\H}{H_0},\;\;\rho_\star=\frac{\rho}{\rho_0},\;\;
\theta_\star=\frac{\theta}{\theta_0},\nonumber
\end{gather}
where the subscripts $0$ denote the corresponding typical values and
$\star$ denotes dimensionless quantities. For convenience, all the
coefficients are assumed to be constants. Thus, the dimensionless form of the system
\eqref{naa}--\eqref{nac} and \eqref{nagg} is obtained by a direct computation:
\begin{align}
&\frac{\partial\rho_\star}{\partial t_\star} +\dv_\star(\rho_\star\u_\star)=0, \nonumber \\
&\rho_\star\frac{d \u_\star}{dt_\star}+\frac{1}{
M^2}\nabla_\star(\rho_\star\theta_\star)
  =C(\nabla_\star \times \H_\star)\times \H_\star+\frac{1}{R}\dv_\star\Psi_\star, \nonumber \\
&\rho_\star\frac{d
\theta_\star}{dt_\star}+(\gamma-1)\rho_\star\theta_\star\dv_\star\u_\star=
\frac{ (\gamma-1)}{R_m}CM^2|\nabla_\star\times\H_\star|^2\nonumber\\
&\qquad \quad
+\frac{ (\gamma-1)M^2}{R} \Psi_\star:\nabla_\star\u_\star+\frac{\gamma}{RP_r}
\Delta_\star \theta_\star,
  \nonumber \\
&\frac{\partial\H_\star}{\partial t_\star}
-\nabla_\star\times(\u_\star\times\H_\star)=\frac{1}{R_m}\nabla_\star\times
(\nabla_\star\times\H_\star),\quad \dv_\star\H_\star=0,\nonumber
\end{align}
where we have used the material derivative
\begin{align}
\frac{d}{dt_\star}=\frac{\partial}{\partial
t_\star}+\u_\star\cdot\nabla_\star,\nonumber
\end{align}
and the new viscous stress tensor
\begin{equation*}
\Psi_\star=2
\mathbb{D}_\star(\u_\star)+\frac{\lambda}{\mu}\dv_\star\u_\star
\;\mathbf{I}_3
\end{equation*}
with
$\mathbb{D}_\star(\u_\star)= (\nabla_\star\u_\star+\nabla_\star\u_\star^\top)/{2}$.

In the above dimensionless system, there are following
dimensionless characteristic parameters:
\begin{eqnarray*} &&
 \mbox{Reynolds number: }
R=\frac{\rho_0u_0L_0}{\mu},\qquad\mbox{Mach number: } M=\frac{u_0}{a_0},   \\
&& \mbox{Prandtl number: } P_r=\frac{c_p\mu}{\kappa},\qquad\mbox{magnetic
Reynolds number: } R_m=\frac{v_0L_0}{\nu},   \\
&&  \mbox{Cowling number: } C=\frac{\mu H_0^2/4\pi\rho_0}{u_0^2},
\end{eqnarray*}
where $c_p$ is the specific heat at constant pressure and
$a_0=\sqrt{ \mathfrak{R}\theta_0}$ is the sound speed. Note
that $\mathfrak{R}=c_p-c_V$ and $\gamma={c_p}/{c_V}$.

\bigskip \noindent
{\bf Acknowledgements:}  The authors are very grateful to the referees for their helpful suggestions. This work was partially done when Fucai Li was
visiting the Institute of Applied Physics and Computational
Mathematics in Beijing. He would like to thank the institute for hospitality.
Jiang was supported by NSFC (Grant No. 40890154) and the National Basic Research Program under the Grant 2011CB309705. Ju was supported by
NSFC (Grant No. 40890154 and 11171035). Li was supported by NSFC (Grant No. 10971094),
PAPD, NCET, and the Fundamental Research Funds for the Central Universities.



\end{document}